\documentclass[DIV=12]{scrartcl}

\usepackage{amsmath}
  \numberwithin{equation}{section}
\usepackage{amssymb}
\usepackage{amsthm}

\theoremstyle{remark}

\usepackage[main=english,german]{babel}
\usepackage[noisbn,noissn,fixlanguage]{babelbib}
 \setbtxfallbacklanguage{english}
\usepackage{cite}%
\usepackage[T1]{fontenc}
\usepackage[final]{hyperref}
\usepackage{lmodern}
\usepackage{mathscinet}
 
\usepackage{mathrsfs}
\usepackage{mathtools}
 \mathtoolsset{mathic=true}
\usepackage{microtype}
\usepackage{pmat}

\usepackage{204arxiv}
\allowdisplaybreaks[1]
\title{\(2\times2\) block representations of the Moore--Penrose inverse and orthogonal projection matrices}
\author{Bernd Fritzsche \and Conrad M\"adler}

\begin{document}
\maketitle

\begin{abstract}
 In this paper, new \tbr{s} of Moore--Penrose inverses for arbitrary complex \taaa{2}{2}{block} matrices are given.
 The approach is based on \tbr{s} of orthogonal projection matrices.
\end{abstract}

\begin{description}
 \item[Keywords] Moore--Penrose inverse, generalized inverses of matrices, block representations, orthogonal projection matrices
 \item[Mathematics Subject Classification (2010)] 15A09 (15A23)
\end{description}

\section{Introduction}
 The aim of this paper is the following:
 Given an arbitrary complex \taaa{\rk{p+q}}{\rk{s+t}}{block} matrix
\beql{Eblock}
 \bE
 =\Mat{a&b\\c&d}
\eeq 
 with \taaa{p}{s}{block} \(a\), we give new \tbr{s} \(\bE^\mpi=\tmat{\alpha&\beta\\\gamma&\delta}\) with \taaa{s}{p}{block} \(\alpha\) of the Moore--Penrose inverse \(\bE^\mpi\) of \(\bE\) as well as new \tbr{s} \(\OPu{\ran{\bE}}=\tmat{e_{11}&e_{12}\\e_{21}&e_{22}}\) with \taaa{p}{p}{block} \(e_{11}\) of the orthogonal projection matrix \(\OPu{\ran{\bE}}\) onto the column space \(\ran{\bE}\) of \(\bE\).
 The block entries should be given by expressions involving the blocks \(a\), \(b\), \(c\), and \(d\) of \(\bE\), where generalized inverses are build only of matrices of the block sizes, \tie{}, with number of rows and columns from the set \(\set{p,q,s,t}\).
 We will see that this goal can be realized (without additional assumptions) by computation of \(\set{1}\)\nobreakdash-inverses of four (\tnnH{}) matrices of the block sizes.

 To the best of our knowledge, except the above mentioned authors Hung/Markham \cite{MR369383} only Miao \cite{MR1102137}, Gro\ss{} \cite{MR1799987}, and Yan \cite{MR3223886} describe explicit \tbr{s} of the Moore--Penrose inverse of arbitrary complex \taaa{2}{2}{partitioned} matrices without making additional assumptions.
 In Miao \cite{MR1102137}, a certain weighted Moore--Penrose inverse is used in the formulas for the block entries of \(\bE^\mpi\).
 In \cite {MR1799987}, Gro\ss{} gives a block representation of the Moore--Penrose inverse of \tnnH{} \taaa{2}{2}{block} matrices.
 Thus, the well-known formulas \(\bE^\mpi=\bE^\ad\rk{\bE\bE^\ad}^\mpi\) and \(\bE^\mpi=\rk{\bE^\ad\bE}^\mpi\bE^\ad\) can then  easily be used to derive a block representation for the Moore--Penrose inverse \(\bE^\mpi\) of an arbitrary block matrix \(\bE\).
 In Yan \cite{MR3223886}, a full rank factorization of \(\bE\) derived from full rank factorizations of the block entries \(a,b,c,d\) is utilized to obtain a \tbr{} of \(\bE^\mpi\).
 
 Assuming certain additional conditions, \teg{}, on column spaces or ranks, several authors derived \tbr{s} of the Moore--Penrose inverse of matrices, see, \teg{} \cite{MR3314341,MR401798,MR3853483}.
 The existence of a Banachiewicz--Schur form for \(\bE^\mpi\) is studied \teg{}\ in \cite{MR1927646,MR2490704}.
 Furthermore, special classes of matrices were considered in this context, see, \teg{}\ \cite{MR2599074,MR469933} for so-called block \(k\)\nobreakdash-circulant matrices.
 Block representations of \(\bE^\mpi\) involving regular transformations, \teg{}, permutations, have also been considered, see \teg{}\ \cite{MR2296072,MR1488625}.
 A representation of the Moore--Penrose inverse of a block column or block row in terms of the block entries is given in \cite{MR2290682}.
 Under a certain rank additivity condition, a \tbr{} of \taaa{m}{n}{partitioned} matrices can be found in \cite{MR1657194} as well.
 Several results on \tbr{s} of partitioned operators are obtained, \teg{}, in \cite{MR2641519,MR2572505,MR2993360,MR2517847,MR2416599}.
 The here mentioned list of references on this topic is not exhaustive.

\section{Notation}
 Throughout this paper let \(m,n,p,q,s,t\) be positive integers.
 We denote by \(\Coo{m}{n}\) the set of all complex \taaa{m}{n}{matrices} and by \(\Co{n}\defeq\Coo{n}{1}\) the set of all column vectors with \(n\) complex entries.
 We write \(\Ouu{m}{n}\) for the zero matrix in \(\Coo{m}{n}\) and \(\Iu{n}\) for the identity matrix in \(\Coo{n}{n}\).
 Let \(\cU\) and \(\cV\) be linear subspaces of \(\Co{n}\).
 If \(\cU\cap\cV=\set{\Ouu{n}{1}}\), then we write \(\cU\oplus\cV\) for the direct sum of \(\cU\) and \(\cV\).
 Let \(\cU^\bot\) be the orthogonal complement of \(\cU\).
 If \(\cU\subseteq\cV\), we use the notation \(\cV\ominus\cU\defeq\cV\cap\rk{\cU^\bot}\).
 We write \(\ran{M}\) and \(\nul{M}\) for the column space and the null space of a complex matrix \(M\).
 Let \(M^\ad\) be the conjugate transpose of a complex matrix \(M\).
 If \(M\) is an arbitrary complex \taaa{m}{n}{matrix}, then there exists a unique complex \taaa{n}{m}{matrix} \(X\) such that the four equations
\begin{align}\label{mpi}
 \text{(1) }MXM&=M,&
 \text{(2) }XMX&=X,&
 \text{(3) }\rk{MX}^\ad&=MX,&
 \text{(4) }\rk{XM}^\ad&=XM
\end{align}
 are fulfilled (see \cite{MR69793}).
 This matrix \(X\) is called the Moore--Penrose inverse of \(M\) and is designated usually by the notation \(M^\mpi\).
 Following~\zitaa{MR1987382}{\cch{1}, \csec{1}, \cdefn{1}{40}}, for each \(M\in\Coo{m}{n}\) we denote by \(M\set{j,k,\dotsc,\ell}\) the set of all \(X\in\Coo{n}{m}\) which satisfy equations \((j),(k),\dotsc,(\ell)\) from the equations (1)--(4) in \eqref{mpi}.
 Each matrix belonging to \(M\set{j,k,\dotsc,\ell}\) is said to be a \emph{\(\set{j,k,\dotsc,\ell}\)\nobreakdash-inverse of \(M\)}.
 
\breml{R1111}
 Let \(\cU\) be a linear subspace of \(\Co{n}\).
 Then there exists a unique complex \taaa{n}{n}{matrix} \(\OPu{\cU}\) such that \(\OPu{\cU}x\in\cU\) and \(x-\OPu{\cU}x\in\cU^\bot\) for all \(x\in\Co{n}\).
 This matrix \(\OPu{\cU}\) is called the orthogonal projection matrix onto \(\cU\).
 If \(P\in\Coo{n}{n}\), then \(P=\OPu{\cU}\) if and only if the three conditions \(P^2=P\) and \(P^\ad=P\) as well as \(\ran{P}=\cU\) are fulfilled.
 Furthermore, the equation \(\OPu{\cU^\bot}=\Iu{n}-\OPu{\cU}\) holds true.
\erem

 Our strategy to give a \tbr{} of the Moore--Penrose inverse \(\bE^\mpi\) of the block matrix \(\bE\) given in \eqref{Eblock} consists of three elementary steps:

\subparagraph{Step~(I)}
 We consider the following factorization problem for orthogonal projection matrices:
 Find a complex \taaa{\rk{s+t}}{\rk{p+q}}{matrix} \(\bR =\tmat{r_{11}&r_{12}\\r_{21}&r_{22}}\) fulfilling \(\OPu{\ran{\bE}}=\bE \bR \) with block entries \(r_{11}\in\Coo{s}{p}\), \(r_{12}\in\Coo{s}{q}\), \(r_{21}\in\Coo{t}{p}\), and \(r_{22}\in\Coo{t}{q}\) expressible explicitly only using the block entries \(a\), \(b\), \(c\), and \(d\) of \(\bE\).
 
\breml{R0811}
 Let \(M\in\Coo{m}{n}\) and let \(X\in\Coo{n}{m}\).
 In view of \rrem{R1111}, then:
\benui
 \il{R0811.a} \(\OPu{\ran{M}}=MX\) if and only if \(X\in M\set{1,3}\).
 \il{R0811.b} \(\OPu{\ran{M^\ad}}=XM\) if and only if \(X\in M\set{1,4}\).
\eenui
\erem

 We constructing a suitable \(\set{1,3}\)\nobreakdash-inverse \(\bR \) of \(\bE\) using:
 
\bthmnl{Urquhart \cite{MR227186}, see also, \teg{}~\zitaa{MR347845}{\cch{1}, \csec{5}, \cthm{3}{47}}}{T1349}
 Let \(M\in\Coo{m}{n}\).
\benui
 \il{T1349.a} Let \(G\defeq MM^\ad\) and let \(G^\gia\in G\set{1}\).
 Then \(M^\ad G^\gia\in M\set{1,2,4}\).
 \il{T1349.b} Let \(H\defeq M^\ad M\) and let \(H^\gia\in H\set{1}\).
 Then \(H^\gia M^\ad\in M\set{1,2,3}\).
\eenui 
\ethm
 
 Applying \rthm{T1349}, we will get an explicit \tbr{} of \(\OPu{\ran{\bE}}=\bE \bR \) in terms of \(a\), \(b\), \(c\), and \(d\).

\subparagraph{Step~(II)}
 Analogous to Step~(I) we construct a suitable complex \taaa{\rk{s+t}}{\rk{p+q}}{matrix} \(\bL=\tmat{\ell_{11}&\ell_{12}\\\ell_{21}&\ell{22}}\) fulfilling \(\bL\in\bE\set{1,4}\) and hence \(\OPu{\ran{\bE^\ad}}=\bL\bE\).

\subparagraph{Step~(III)}  
 With the matrices \(\bL\) and \(\bR \) we apply:
 
\bthmnl{Urquhart \cite{MR227186}, see \teg{}~\zitaa{MR347845}{\cch{1}, \csec{5}, \cthm{4}{48}}}{T1411}
 If \(M\in\Coo{m}{n}\), then \(M^\giad MM^\giac =M^\mpi\) for every choice of \(M^\giac \in M\set{1,3}\) and \(M^\giad \in M\set{1,4}\).
\ethm

 Regarding \rrem{R0811}, \rthm{T1411} admits the following reformulation:
\breml{R1509}
 Let \(M\in\Coo{m}{n}\) and let \(L\in\Coo{n}{m}\) and \(R\in\Coo{n}{m}\) be such that \(LM=\OPu{\ran{M^\ad}}\) and \(MR=\OPu{\ran{M}}\).
 Then \(M^\mpi=LMR\).
\erem

 Consider an additive decomposition \(M=U+V\) of an \taaa{m}{n}{matrix} \(M\) with two \taaa{m}{n}{matrices} \(U\) and \(V\), fulfilling \(UV^\ad=\Ouu{m}{m}\).
 In this situation, a result of Cline \cite{MR180572} is applicable to obtain a non-trivial representation of \(M^\mpi\) as a sum of \(U^\mpi\) and a further matrix.
 By Hung/Markham \cite{MR369383} a decomposition \(\bE=\bU+\bV\) with \(\bU=\tmat{a&0\\c&0}\) and \(\bV=\tmat{0&c\\0&d}\) is used in this way to derive a \tbr{} of \(\bE^\mpi\) involving only Moore--Penrose inverses of block size matrices.
 
 Regarding \rremss{R0811}{R1111} as well as \eqref{mpi}, the orthogonal projection matrix \(Q\defeq\OPu{\ran{U^\ad}}\) fulfills \(UQ=UU^\mpi U=U\) and \(VQ=\rk{QV^\ad}^\ad=\rk{U^\mpi UV^\ad}^\ad=\Ouu{m}{n}\).
 Consequently, \(MQ=U\) and \(M{\rk{\Iu{n}-Q}}=V\).
 Conversely, given \(M\in\Coo{m}{n}\) and an orthogonal projection matrix \(Q\in\Coo{n}{n}\), then it is readily checked that \(U\defeq MQ\) and \(V\defeq M\rk{\Iu{n}-Q}\) fulfill \(M=U+V\) and \(UV^\ad=\Ouu{m}{m}\).
 Thus, every decomposition \(M=U+V\) with \(UV^\ad=\Ouu{m}{m}\) can be written as \(M=MQ+M\rk{\Iu{n}-Q}\) with some orthogonal projection matrix \(Q\in\Coo{n}{n}\) occurring on the right-hand side and vice versa.
 Although not explicitly using Cline's theorem, our investigations involve an analogous decomposition, namely \(\bE=\OPu{\cS}\bE+\rk{\Iu{p+q}-\OPu{\cS}}\bE\) with the orthogonal projection matrix \(\OPu{\cS}\) onto the linear subspace \(\cS\) spanned by the first \(s\) columns of \(\bE\) occurring on the left-hand side (see \rlem{L1354}).

\section{Main results}
 We consider a complex \taaa{\rk{p+q}}{\rk{s+t}}{matrix} \(\bE\).
 Let \eqref{Eblock} be the \tbr{} of \(\bE\) with \taaa{p}{s}{block} \(a\).
 Setting
\begin{align}\label{STYZ}
 Y&\defeq\mat{a,b},&
 Z&\defeq\mat{c,d},&
 S&\defeq\Mat{a\\c},&
 T&\defeq\Mat{b\\d}
\end{align}
 then
\begin{align}\label{ESZblock}
 \bE&=\Mat{Y\\Z},&
 \bE&=\mat{S,T}.
\end{align}
 Let
\begin{align}
 \mu&\defeq aa^\ad+bb^\ad,&
 \sigma&\defeq a^\ad a+c^\ad c,\label{syab}\\
 \zeta&\defeq cc^\ad+dd^\ad,&
 \tau&\defeq b^\ad b+d^\ad d,\notag\\
 \rho&\defeq ca^\ad+db^\ad,&
 \lambda&\defeq a^\ad b+c^\ad d.\label{trab}
\end{align}
 In view of \eqref{STYZ}, then
\begin{align}
 \mu&=YY^\ad,&
 \sigma&=S^\ad S,\label{sy}\\
 \zeta&=ZZ^\ad,&
 \tau&=T^\ad T,\label{zt}\\
 \rho&=ZY^\ad,&
 \lambda&=S^\ad T.\label{tr}
\end{align}
 Choose \(\mu^\gia\in\mu\set{1}\) and \(\sigma^\gia\in\sigma\set{1}\).
 Regarding \eqref{sy}, then \rthm{T1349} shows that
\begin{align}
 Y^\giabd&\defeq Y^\ad\mu^\gia,&
 S^\giabc&\defeq\sigma^\gia S^\ad\label{bbYbbS}
\intertext{fulfill}
 Y^\giabd&\in Y\set{1,2,4},&
 S^\giabc&\in S\set{1,2,3}.\label{YS123}
\end{align}
 Let
\begin{align}
 \phi&\defeq c-\rk{ca^\ad+db^\ad}\mu^\gia a,&
 \psi&\defeq d-\rk{ca^\ad+db^\ad}\mu^\gia b,\label{fg}\\
 \eta&\defeq b-a\sigma^\gia\rk{a^\ad b+c^\ad d},&
 \theta&\defeq d-c\sigma^\gia\rk{a^\ad b+c^\ad d}.\label{et}
\end{align}
 Because of \eqref{bbYbbS}, \eqref{tr}, \eqref{STYZ}, \eqref{trab}, \eqref{fg}, and \eqref{et}, then
\begin{align}\label{VWblock}
 V&\defeq\mat{\phi,\psi}&
&\text{and}&
 W&\defeq\Mat{\eta\\\theta}
\end{align}
 admit the representations
\beql{V}\begin{split}
 Z\rk{\Iu{s+t}-Y^\giabd Y}
 &=Z-ZY^\ad\mu^\gia Y
 =Z-\rho\mu^\gia Y\\
 &=\mat{c,d}-\rho\mu^\gia\mat{a,b}
 =\mat{\phi,\psi}
 =V
\end{split}\eeq
 and
\beql{W}\begin{split}
 \rk{\Iu{p+q}-SS^\giabc}T
 &=T-S\sigma^\gia S^\ad T
 =T-S\sigma^\gia\lambda\\
 &=\Mat{b\\d}-\Mat{a\\c}\sigma^\gia\lambda
 =\Mat{\eta\\\theta}
 =W.
\end{split}\eeq
 Using \eqref{V}, \eqref{W}, \eqref{YS123}, \rremss{R0811}{R1111}, and \(\ek{\ran{Y^\ad}}^\bot=\nul{Y}\), we can infer
\begin{align}\label{VWOP}
 V&=Z\OPu{\nul{Y}},&
 W&=\OPu{\ek{\ran{S}}^\bot}T.
\end{align}
 Let
\begin{align}
 \nu&\defeq\phi\phi^\ad+\psi\psi^\ad,&
 \omega&\defeq\eta^\ad\eta+\theta^\ad\theta.\label{noab}
\end{align}
 In view of \eqref{VWblock}, \eqref{VWOP}, and \rrem{R1111} then
\begin{align}
 \nu&=VV^\ad=Z\OPu{\nul{Y}}Z^\ad=VZ^\ad,&
 \omega&=W^\ad W=T^\ad\OPu{\ek{\ran{S}}^\bot}T=T^\ad W.\label{no}
\end{align}
 Choose \(\nu^\gia\in\nu\set{1}\) and \(\omega^\gia\in\omega\set{1}\).
 Regarding \eqref{no}, then \rthm{T1349} shows that
\begin{align}
 V^\giabd&\defeq V^\ad\nu^\gia&
&\text{and}&
 W^\giabc&\defeq\omega^\gia W^\ad\label{bbVbbW}
\intertext{fulfill}
 V^\giabd&\in V\set{1,2,4}&
&\text{and}&
 W^\giabc&\in W\set{1,2,3}.\label{VW123}
\end{align}
 Obviously, we have \(\mu\in\Cggo{p}\), \(\sigma\in\Cggo{s}\), \(\zeta\in\Cggo{q}\), \(\tau\in\Cggo{t}\), \(\nu\in\Cggo{q}\), \(\omega\in\Cggo{t}\), \(\rho\in\Coo{q}{p}\), and \(\lambda\in\Coo{s}{t}\), where \(\Cggo{n}\) denotes the set of all \tnnH{} complex \taaa{n}{n}{matrices}.

\breml{R1005}
 Let
\begin{align}\label{LR}
 \bL&\defeq\mat*{\rk{\Iu{s+t}-V^\giabd Z}Y^\giabd,V^\giabd},&
 \bR &\defeq\Mat{S^\giabc\rk{\Iu{p+q}-TW^\giabc}\\W^\giabc}.
\end{align}
 Regarding \eqref{LR}, \eqref{bbVbbW}, \eqref{bbYbbS}, \eqref{tr}, \eqref{STYZ}, and \eqref{VWblock}, then
\[\begin{split}
 \bL
 &=\mat*{\rk{\Iu{s+t}-V^\ad\nu^\gia Z}Y^\ad\mu^\gia,V^\ad\nu^\gia}
 =\mat*{\rk{Y^\ad-V^\ad\nu^\gia ZY^\ad}\mu^\gia,V^\ad\nu^\gia}\\
 &=\mat*{\rk{Y^\ad-V^\ad\nu^\gia\rho}\mu^\gia,V^\ad\nu^\gia}
 =\mat{Y^\ad\mu^\gia,\Ouu{\rk{s+t}}{q}}+\mat{-V^\ad\nu^\gia\rho\mu^\gia,V^\ad\nu^\gia}\\
 &=Y^\ad\mu^\gia\mat{\Iu{p},\Ouu{p}{q}}+V^\ad\nu^\gia\mat{-\rho\mu^\gia,\Iu{q}}\\
 &=\Mat{a^\ad\\b^\ad}\mu^\gia\mat{\Iu{p},\Ouu{p}{q}}+\Mat{\phi^\ad\\\psi^\ad}\nu^\gia\mat{-\rho\mu^\gia,\Iu{q}}
 =\Mat{
 \ell_{11}&\ell_{12}\\
 \ell_{21}&\ell_{22}},
\end{split}\]
 where
\begin{align}\label{l}
 \ell_{12}&\defeq\phi^\ad\nu^\gia,&
 \ell_{11}&\defeq\rk{a^\ad-\ell_{12}\rho}\mu^\gia,&
 \ell_{22}&\defeq\psi^\ad\nu^\gia,&
 \ell_{21}&\defeq\rk{b^\ad-\ell_{22}\rho}\mu^\gia,
\end{align}
 and
\[\begin{split}
 \bR 
 &=\Mat{\sigma^\gia S^\ad\rk{\Iu{p+q}-T\omega^\gia W^\ad}\\\omega^\gia W^\ad}
 =\Mat{\sigma^\gia\rk{S^\ad-S^\ad T\omega^\gia W^\ad}\\
\omega^\gia W^\ad}\\
 &=\Mat{\sigma^\gia\rk{S^\ad-\lambda\omega^\gia W^\ad}\\
\omega^\gia W^\ad}
 =\Mat{\sigma^\gia S^\ad\\
\Ouu{t}{\rk{p+q}}}+\Mat{-\sigma^\gia\lambda\omega^\gia W^\ad\\\omega^\gia W^\ad}\\
 &=\Mat{\Iu{s}\\\Ouu{t}{s}}\sigma^\gia S^\ad+\Mat{-\sigma^\gia\lambda\\
\Iu{t}}\omega^\gia W^\ad\\
 &=\Mat{\Iu{s}\\\Ouu{t}{s}}\sigma^\gia\mat{a^\ad,c^\ad}+\Mat{-\sigma^\gia\lambda\\
\Iu{t}}\omega^\gia\mat{\eta^\ad,\theta^\ad}
 =\Mat{
 r_{11}&r_{12}\\
 r_{21}&r_{22}},
\end{split}\]
 where
\begin{align}\label{r}
 r_{21}&\defeq\omega^\gia\eta^\ad,&
 r_{11}&\defeq\sigma^\gia\rk{a^\ad-\lambda r_{21}},&
 r_{22}&\defeq\omega^\gia\theta^\ad,&
 r_{12}&\defeq\sigma^\gia\rk{c^\ad-\lambda r_{22}}.
\end{align}
\erem

\bleml{L1354}
 Let \(\cE\defeq\ran{\bE}\), let \(\cS\defeq\ran{S}\), and let \(\cW\defeq\ran{W}\).
 Then \(\cW=\cE\ominus\cS\) and \(\OPu{\cE}=\OPu{\cS}+\OPu{\cW}=\bE \bR \).
\elem
\bproof
 We first check \(\cW=\cE\cap\rk{\cS^\bot}\).
 Because of \eqref{VW123}, \eqref{YS123}, and \rremp{R0811}{R0811.a}, we have \(\OPu{\cW}=WW^\giabc\) and \(\OPu{\cS}=SS^\giabc\).
 According to \rrem{R1111}, hence \(\Iu{p+q}-SS^\giabc=\OPu{\cS^\bot}\).
 By virtue of \eqref{W}, then \(W=\OPu{\cS^\bot}T\) follows.
 From \rrem{R1111} we know \(\ran{\OPu{\cS^\bot}}=\cS^\bot\).
 Consequently, \(\cW\subseteq\cS^\bot\).
 Regarding \eqref{ESZblock} and \eqref{W}, we have \(\cS\subseteq\cE\) and
\[
 W
 =\rk{\Iu{p+q}-SS^\giabc}T
 =T-SS^\giabc T
 =\mat{S,T}\Mat{-S^\giabc T\\\Iu{t}}
 =\bE\Mat{-S^\giabc T\\\Iu{t}},
\]
 implying \(\cW\subseteq\cE\).
 Thus, \(\cW\subseteq\cE\cap\rk{\cS^\bot}\) is proved.
 Now we consider an arbitrary \(\bw\in\cE\cap\rk{\cS^\bot}\).
 Then \(\bw\in\cE\); so there exists some \(\bv\in\Co{s+t}\) with \(\bw=\bE\bv\).
 Let \(\bv=\tmat{x\\y}\) be the \tbr{} of \(\bv\) with \(x\in\Co{s}\).
 Regarding \eqref{ESZblock}, then \(\bw=Sx+Ty\).
 In view of \(\bw\in\cS^\bot\) and \(Sx\in\cS\), furthermore \(\OPu{\cS^\bot}\bw=\bw\) and \(\OPu{\cS^\bot}Sx=\Ouu{\rk{p+q}}{1}\).
 Taking additionally into account \(W=\OPu{\cS^\bot}T\), we obtain then
\[
 \bw
 =\OPu{\cS^\bot}\bw
 =\OPu{\cS^\bot}\rk{Sx+Ty}
 =\OPu{\cS^\bot}Ty
 =Wy,
\]
 implying \(\bw\in\cW\).
 Thus, we have also shown \(\cE\cap\rk{\cS^\bot}\subseteq\cW\).
 Therefore, \(\cW=\cE\cap\rk{\cS^\bot}\) holds true.
 Since \(\cS\subseteq\cE\), hence \(\cW=\cE\ominus\cS\).
 Consequently, \(\OPu{\cW}=\OPu{\cE}-\OPu{\cS}\) follows (see, \teg{}~\zitaa{MR566954}{\cthm{4.30(c)}{82}}).
 Thus, \(\OPu{\cE}=\OPu{\cS}+\OPu{\cW}\).
 Taking additionally into account \(\OPu{\cS}=SS^\giabc\) and \(\OPu{\cW}=WW^\giabc\) as well as \eqref{W}, \eqref{ESZblock}, and \eqref{LR}, then we can conclude
\[\begin{split}
 \OPu{\cE}
 &=SS^\giabc+WW^\giabc
 =SS^\giabc+\rk{\Iu{p+q}-SS^\giabc}TW^\giabc\\
 &=SS^\giabc\rk{\Iu{p+q}-TW^\giabc}+TW^\giabc
 =\mat{S,T}\Mat{S^\giabc\rk{\Iu{p+q}-TW^\giabc}\\W^\giabc}
 =\bE \bR.\qedhere
\end{split}\]
\eproof

 The following result can be proved analogously.
 We omit the details.

\bleml{L0921}
 Let \(\tilde{\cE}\defeq\ran{\bE^\ad}\), let \(\tilde{\cY}\defeq\ran{Y^\ad}\), and let \(\tilde{\cV}\defeq\ran{V^\ad}\).
 Then \(\tilde{\cV}=\tilde{\cE}\ominus\tilde{\cY}\) and \(\OPu{\tilde{\cE}}=\OPu{\tilde{\cY}}+\OPu{\tilde{\cV}}=\bL\bE\).
\elem

\breml{R0944}
 From \rlem{L1354} and \rremp{R0811}{R0811.a} we can infer \(R\in\bE\set{1,3}\), whereas \rlem{L0921} and \rremp{R0811}{R0811.b} yield \(L\in\bE\set{1,4}\).
\erem

 Now we obtain the announced \tbr{s} of orthogonal projection matrices.

\bpropl{P1819n}
 Let \(\bE\) be a complex \taaa{\rk{p+q}}{\rk{s+t}}{matrix} and let \eqref{Eblock} be the \tbr{} of \(\bE\) with \taaa{p}{s}{block} \(a\).
 Let \(S\) be given by \eqref{STYZ}.
 Let \(\sigma\) and \(\lambda\) be given by \eqref{syab} and \eqref{trab}.
 Let \(\sigma^\gia\in\sigma\set{1}\).
 Let \(\eta,\theta\) and \(W\) be given by \eqref{et} and \eqref{VWblock}.
 Let \(\omega\) be given by \eqref{noab} and let \(\omega^\gia\in\omega\set{1}\).
 Then
\begin{align*}
 \OPu{\ran{\bE}}
 &=S\sigma^\gia S^\ad+W\omega^\gia W^\ad
 =
 \begin{pmat}[{|}]
  a\sigma^\gia a^\ad+\eta\omega^\gia\eta^\ad&a\sigma^\gia c^\ad+\eta\omega^\gia\theta^\ad\cr\-
  c\sigma^\gia a^\ad+\theta\omega^\gia\eta^\ad&c\sigma^\gia c^\ad+\theta\omega^\gia\theta^\ad\cr
 \end{pmat}.
\end{align*} 
\eprop
\bproof
 In the proof of \rlem{L1354}, we have already shown \(\OPu{\ran{\bE}}=SS^\giabc+WW^\giabc\).
 Taking additionally into account \eqref{bbYbbS}, \eqref{bbVbbW}, \eqref{STYZ}, and \eqref{VWblock}, the assertions follow.
\eproof

 Now we are able to prove a \(2\times2\) \tbr{} of the Moore--Penrose inverse.

\bthml{T1009}
 Let \(\bE\) be a complex \taaa{\rk{p+q}}{\rk{s+t}}{matrix} and let \eqref{Eblock} be the \tbr{} of \(\bE\) with \taaa{p}{s}{block} \(a\).
 Let \(Y,S\) be given by \eqref{STYZ}.
 Let \(\mu,\sigma\) and \(\rho,\lambda\) be given by \eqref{syab} and \eqref{trab}.
 Let \(\mu^\gia\in\mu\set{1}\) and let \(\sigma^\gia\in\sigma\set{1}\).
 Let \(\phi,\psi\) and \(\eta,\theta\) be given by \eqref{fg} and \eqref{et}.
 Let \(V\) and \(W\) be given by \eqref{VWblock}.
 Let \(\nu\) and \(\omega\) be given by \eqref{noab}.
 Let \(\nu^\gia\in\nu\set{1}\) and let \(\omega^\gia\in\omega\set{1}\).
 Then
\begin{multline*}
 \bE^\mpi
 =\bL\bE \bR 
 =\Mat{\alpha&\beta\\\gamma&\delta}
 =\rk{Y^\ad-V^\ad\nu^\gia\rho}\mu^\gia a\sigma^\gia\rk{S^\ad-\lambda\omega^\gia W^\ad}+\rk{Y^\ad-V^\ad\nu^\gia\rho}\mu^\gia b\omega^\gia W^\ad\\
+V^\ad\nu^\gia c\sigma^\gia\rk{S^\ad-\lambda\omega^\gia W^\ad}+V^\ad\nu^\gia d\omega^\gia W^\ad
\end{multline*}
 with
\begin{align*}
 \alpha&\defeq\ell_{11}ar_{11}+\ell_{11}br_{21}+\ell_{12}cr_{11}+\ell_{12}dr_{21}&
 \beta&\defeq\ell_{11}ar_{12}+\ell_{11}br_{22}+\ell_{12}cr_{12}+\ell_{12}dr_{22}\\
 \gamma&\defeq\ell_{21}ar_{11}+\ell_{21}br_{21}+\ell_{22}cr_{11}+\ell_{22}dr_{21}&
 \delta&\defeq\ell_{21}ar_{12}+\ell_{21}br_{22}+\ell_{22}cr_{12}+\ell_{22}dr_{22}
\end{align*}
 where, for each \(j,k\in\set{1,2}\), the matrices \(\ell_{jk}\) and \(r_{jk}\) are given by \eqref{l} and \eqref{r}, \tresp{}
\ethm
\bproof
 According to \rlemss{L0921}{L1354}, we have \(\OPu{\ran{\bE^\ad}}=\bL\bE\) and \(\OPu{\ran{\bE}}=\bE \bR \).
 Thus, we can apply \rrem{R1509} to obtain \(\bE^\mpi=\bL\bE \bR \).
 Using \rrem{R1005} and \eqref{Eblock}, we furthermore obtain
\begin{multline*}
 \bL\bE \bR 
 =\mat*{\rk{Y^\ad-V^\ad\nu^\gia\rho}\mu^\gia,V^\ad\nu^\gia}
 \Mat{a&b\\c&d}
 \Mat{\sigma^\gia\rk{S^\ad-\lambda\omega^\gia W^\ad}\\
\omega^\gia W^\ad}\\
 =\rk{Y^\ad-V^\ad\nu^\gia\rho}\mu^\gia a\sigma^\gia\rk{S^\ad-\lambda\omega^\gia W^\ad}+\rk{Y^\ad-V^\ad\nu^\gia\rho}\mu^\gia b\omega^\gia W^\ad\\
+V^\ad\nu^\gia c\sigma^\gia\rk{S^\ad-\lambda\omega^\gia W^\ad}+V^\ad\nu^\gia d\omega^\gia W^\ad
\end{multline*}
 as well as
\[\begin{split}
 \bL\bE \bR 
 &=
 \Mat{
 \ell_{11}&\ell_{12}\\
 \ell_{21}&\ell_{22}}
 \Mat{
 a&b\\
 c&d
 }
 \Mat{
 r_{11}&r_{12}\\
 r_{21}&r_{22}
 }
 =\Mat{\alpha&\beta\\\gamma&\delta}.\qedhere
\end{split}\]
\eproof

\section{Examples for consequences}
 In this section, we give some examples of applications of the \tbr{s} of orthogonal projection matrices and Moore--Penrose inverses given in \rprop{P1819n} and \rthm{T1009}.
 In order to avoid lengthy formulas, we give only hints for computations.
 
\bexal{R13}
 Let \(N\) be a positive integer and let \(\cU\) and \(\cV\) be two complementary linear subspaces of \(\Co{N}\), \tie{}, the subspaces \(\cU\) and \(\cV\) fulfill \(\cU\oplus\cV=\Co{N}\).
 Then there exists a unique complex \taaa{N}{N}{matrix} \(\PPuu{\cU}{\cV}\) such that \(\PPuu{\cU}{\cV}x=u\) for all \(x\in\Co{N}\), where \(x=u+v\) is the unique representation of \(x\) with \(u\in\cU\) and \(v\in\cV\).
 This matrix \(\PPuu{\cU}{\cV}\) is called the oblique projection matrix on \(\cU\) along \(\cV\) and admits the representations
\beql{R13.1}
 \PPuu{\cU}{\cV}
 =\rk{\OPu{\cV^\bot}\OPu{\cU}}^\mpi
 =\ek*{\rk{\Iu{N}-\OPu{\cV}}\OPu{\cU}}^\mpi,
\eeq
 see \cite{MR347845} or also \zitaa{MR1987382}{\cch{2}, \csec{7}, \cex{60}{80}, \cform{(80)}{80}}.
 Assume that \(N=p+q\) and that \(\cU=\ran{\bE}\) and \(\cV=\ran{\bF}\) for two matrices \(\bE\in\Coo{\rk{p+q}}{\rk{s+t}}\) and \(\bF\in\Coo{\rk{p+q}}{\rk{m+n}}\) with \tbr{s} \eqref{Eblock} and \(\bF=\tmat{e&f\\g&h}\), where \(a\) is a \taaa{p}{s}{block} and \(e\) is a \taaa{p}{m}{block}.
 Then \eqref{R13.1} together with \rprop{P1819n} and \rthm{T1009} could be used to obtain a \tbr{} of \(\PPuu{\cU}{\cV}\) in terms of \(a,b,c,d\) and \(e,f,g,h\).
\eexa

\bexal{E1216}
 Because \(\cU\oplus\rk{\cU^\bot}=\Co{n}\) holds true for every linear subspace \(\cU\) of \(\Co{n}\), the Moore--Penrose inverse of matrices is a special case of the uniquely determined \(\set{1,2}\)\nobreakdash-inverse with simultaneously prescribed column space and null space.
 More precisely, if \(M\in\Coo{m}{n}\) and linear subspaces \(\cU\) of \(\Co{m}\) and \(\cV\) of \(\Co{n}\) with \(\ran{M}\oplus\cU=\Co{m} \) and \(\nul{M}\oplus\cV=\Co{n}\) are given, then there exists a unique complex \taaa{n}{m}{matrix} \(X\) such that the four conditions
\begin{align*} 
 MXM&=M,&
 XMX&=X,&
 \ran{X}&=\cV,&
&\text{and}&
 \nul{X}&=\cU
\end{align*}
 are fulfilled.
 This matrix \(X\) is denoted by \(\giuu{M}{\cU}{\cV}\) and admits the representations
\beql{E1216.1}
 \giuu{M}{\cV}{\cU}
 =\PPuu{\cV}{\nul{M}}M^\gia\PPuu{\ran{M}}{\cU}
 =\rk{\OPu{\cV^\bot}\OPu{\nul{M}}}^\mpi M^\gia\rk{\OPu{\ek{\ran{M}}^\bot}\OPu{\cU}}^\mpi
\eeq
 with every \(M^\gia\in M\set{1}\) (see, \teg{} \zitaa{MR1987382}{\cch{2}, \csec{6}, \cthm{12}{71}}).
 In particular, \eqref{E1216.1} is valid for \(M^\gia=M^\mpi\).
 Furthermore, \(M^\mpi=\giuu{M}{\ek{\nul{M}}^\bot}{\ek{\ran{M}}^\bot}\).
 \rexa{R13} could be used to obtain a \tbr{} of \(\giuu{M}{\cV}{\cU}\), if the subspaces \(\cU\) and \(\cV\) are given as column spaces of certain matrices, partitioned accordingly to a \tbr{} of \(M\), and if a matrix \(M^\gia\in M\set{1}\) the corresponding \tbr{} of which is known.
 (In particular, for \(M^\gia=M^\mpi\) one can use \rthm{T1009}.)
\eexa

\section{An alternative approach}
 In this final section, we give alternative representations of the matrices \(\bL\) and \(\bR\) occurring in \rthm{T1009} and \rlemss{L0921}{L1354}.
 Utilizing these representations, further \tbr{s} of the Moore--Penrose inverse \(\bE^\mpi\) could possibly be obtained, in particular, in the case of \(\bE\) satisfying additional conditions.
 We will not pursue this direction any further here.
 We continue to use the notations given above.

\blemnl{Rohde \cite{MR190161}, see, \teg{}~\zitaa{MR347845}{\cch{5}, \csec{2}, \cex{10(a)}{182}}}{P1358}
 Let \(M\in\Cggo{\rk{p+q}}\) and let \(M=\tmat{m_{11}&m_{12}\\m_{21}&m_{22}}\) be the \tbr{} of \(M\) with \taaa{p}{p}{block} \(m_{11}\).
 Let \(m_{11}^\gia\in m_{11}\set{1}\) and let \(\varsigma\defeq m_{22}-m_{21}m_{11}^\gia m_{12}\).
 Let \(\varsigma^\gia\in\varsigma\set{1}\).
 Then
\[
 M^\gia
 \defeq
 \begin{pmat}[{|}]
  m_{11}^\gia+m_{11}^\gia m_{12}\varsigma^\gia m_{21}m_{11}^\gia&-m_{11}^\gia m_{12}\varsigma^\gia\cr\-
  -\varsigma^\gia m_{21}m_{11}^\gia&\varsigma^\gia\cr
 \end{pmat}
\]
 belongs to \(M^\gia\in M\set{1}\).
\elem

\breml{R1406}
 Regarding \eqref{no}, \eqref{V}, \eqref{W}, \eqref{bbYbbS}, \eqref{zt}, and, \eqref{tr}, we can infer
\begin{align}
 \nu
 &=VZ^\ad
 =Z\rk{\Iu{s+t}-Y^\giabd Y}Z^\ad
 =Z\rk{\Iu{s+t}-Y^\ad\mu^\gia Y}Z^\ad
 =\zeta-\rho\mu^\gia \rho^\ad\label{xiZV}
\intertext{and}
 \omega
 &=T^\ad W 
 =T^\ad\rk{\Iu{p+q}-SS^\giabc}T
 =T^\ad\rk{\Iu{p+q}-S\sigma^\gia S^\ad}T
 =\tau-\lambda^\ad\sigma^\gia\lambda.\label{piWT}
\end{align}
\erem

\bleml{C1010}
 Let
\begin{align}\label{GH}
 \bG&\defeq\bE\bE^\ad&
&\text{and}&
 \bH&\defeq\bE^\ad\bE.
\end{align}
 Let \(\mu^\gia\in\mu\set{1}\) and \(\nu^\gia\in\nu\set{1}\) and let \(\sigma^\gia\in\sigma\set{1}\) and \(\omega^\gia\in\omega\set{1}\).
 Then the matrices
\begin{align}
 \bG^\gia
 &\defeq
 \begin{pmat}[{|}]
  \mu^\gia+\mu^\gia \rho^\ad\nu^\gia \rho\mu^\gia&-\mu^\gia \rho^\ad\nu^\gia\cr\-
  -\nu^\gia \rho\mu^\gia&\nu^\gia\cr
 \end{pmat}\label{G1}\\
\intertext{and}
 \bH^\gia
 &\defeq
 \begin{pmat}[{|}]
  \sigma^\gia+\sigma^\gia \lambda\omega^\gia \lambda^\ad\sigma^\gia&-\sigma^\gia \lambda\omega^\gia\cr\-
  -\omega^\gia \lambda^\ad\sigma^\gia&\omega^\gia\cr
 \end{pmat}\label{H1}
\end{align}
 fulfill \(\bG^\gia\in\bG\set{1}\) and \(\bH^\gia\in\bH\set{1}\).
\elem
\bproof
 Clearly, \(\bG\in\Cggo{\rk{p+q}}\) and \(\bH\in\Cggo{\rk{s+t}}\).
 Regarding \eqref{ESZblock}, \eqref{sy}, and \eqref{zt}, we have \(\bG=\tmat{Y\\Z}\mat{Y^\ad,Z^\ad}=\tmat{YY^\ad&YZ^\ad\\ZY^\ad&ZZ^\ad}=\tmat{\mu&\rho^\ad\\\rho&\zeta}\) and \(\bH=\tmat{S^\ad\\T^\ad}\mat{S,T}=\tmat{S^\ad S&S^\ad T\\T^\ad S&T^\ad T}=\tmat{\sigma&\lambda\\\lambda^\ad&\tau}\).
 Taking additionally into account \rrem{R1406}, thus from \rlem{P1358} the assertions immediately follow.
\eproof

 Finally, in the following result, we not only get new representations for the matrices \(\bL\) and \(\bR\) occurring in \rthm{T1009} and \rlemss{L0921}{L1354}, but also obtain their belonging to the set \(\bE\set{1,2,4}\) and \(\bE\set{1,2,3}\), \tresp{}, thereby improving \rrem{R0944}.

\bleml{R1401}
 The matrices \(\bL\) and \(\bR\) admit the representations \(\bL=\bE^\ad\bG^\gia\) and \(\bR=\bH^\gia\bE^\ad\) and fulfill \(\bL\in\bE\set{1,2,4}\) and \(\bR\in\bE\set{1,2,3}\).
\elem
\bproof
 Using \eqref{YS123} and \rremss{R0811}{R1111}, we have \(\rk{Y^\giabd Y}^\ad=\OPu{\ran{Y^\ad}}^\ad=\OPu{\ran{Y^\ad}}=Y^\giabd Y\) and, analogously, \(\rk{SS^\giabc}^\ad=SS^\giabc\).
 Taking additionally into account \eqref{V}, \eqref{W}, \eqref{bbYbbS}, and \eqref{tr}, we thus can conclude
\begin{align}
 V^\ad
 &=\rk{\Iu{s+t}-Y^\giabd Y}Z^\ad
 =Z^\ad-Y^\ad\mu^\gia YZ^\ad
 =Z^\ad-Y^\ad\mu^\gia\rho^\ad\label{V*}
\intertext{and}
 W^\ad
 &=T^\ad\rk{\Iu{p+q}-SS^\giabc}
 =T^\ad-T^\ad S\sigma^\gia S^\ad
 =T^\ad-\lambda^\ad\sigma^\gia S^\ad.\label{W*}
\end{align}
 Regarding \eqref{ESZblock}, \eqref{G1}, \eqref{H1}, \eqref{xiZV}, \eqref{piWT}, \eqref{V*}, \eqref{W*}, and \rrem{R1005}, we obtain
\beql{EG}\begin{split}
 \bE^\ad\bG^\gia
 &=\mat{Y^\ad,Z^\ad}\rk*{\Mat{\Iu{p}\\\Ouu{q}{p}}\mu^\gia\mat{\Iu{p},\Ouu{p}{q}}+\Mat{-\mu^\gia\rho^\ad\\\Iu{q}}\nu^\gia\mat{-\rho\mu^\gia,\Iu{q}}}\\
 &=Y^\ad\mu^\gia\mat{\Iu{p},\Ouu{p}{q}}+\rk{Z^\ad-Y^\ad\mu^\gia\rho^\ad}\nu^\gia\mat{-\rho\mu^\gia,\Iu{q}}
 =\bL
\end{split}\eeq
 and
\beql{HE}\begin{split}
 \bH^\gia\bE^\ad
 &=\rk*{\Mat{\Iu{s}\\\Ouu{t}{s}}\sigma^\gia\mat{\Iu{s},\Ouu{s}{t}}+\Mat{-\sigma^\gia \lambda\\\Iu{t}}\omega^\gia\mat{-\lambda^\ad\sigma^\gia,\Iu{t}}}
 \Mat{S^\ad\\T^\ad}\\
 &=\Mat{\Iu{s}\\\Ouu{t}{s}}\sigma^\gia S^\ad+\Mat{-\sigma^\gia \lambda\\\Iu{t}}\omega^\gia\rk{T^\ad-\lambda^\ad\sigma^\gia S^\ad}
 =\bR.
\end{split}\eeq
 According to \rlem{C1010}, we have \(\bG^\gia\in\bG\set{1}\) and \(\bH^\gia\in\bH\set{1}\).
 Taking additionally into account \eqref{GH}, \eqref{EG}, and \eqref{HE}, then \(\bL\in\bE\set{1,2,4}\) and \(\bR\in\bE\set{1,2,3}\) follow from \rthm{T1349}.
\eproof

 Observe that \rlem{R1401} in connection with \rthm{T1009} yields a factorization \(\bE^\mpi=\bL\bE\bR\) with particular matrices \(\bL\in\bE\set{1,2,4}\) and \(\bR\in\bE\set{1,2,3}\).
 This gives a special factorization of the kind mentioned in Urquhart's result (\rthm {T1411}), whereby all matrices can be expressed explicitly in terms of the block entries \(a,b,c,d\) of the given matrix \(\bE\).

\bibliography{204arxiv}

\begin{thebibliography}{10}
  \providebibliographyfont{name}{}%
  \providebibliographyfont{lastname}{}%
  \providebibliographyfont{title}{\emph}%
  \providebibliographyfont{jtitle}{\btxtitlefont}%
  \providebibliographyfont{etal}{\emph}%
  \providebibliographyfont{journal}{}%
  \providebibliographyfont{volume}{}%
  \providebibliographyfont{ISBN}{\MakeUppercase}%
  \providebibliographyfont{ISSN}{\MakeUppercase}%
  \providebibliographyfont{url}{\url}%
  \providebibliographyfont{numeral}{}%
  \expandafter\btxselectlanguage\expandafter {\btxfallbacklanguage}

\expandafter\btxselectlanguage\expandafter {\btxfallbacklanguage}
\bibitem {MR2290682}
\btxnamefont {\btxlastnamefont {Baksalary},~J.\btxfnamespaceshort K.}
  \btxandshort {.}\ \btxnamefont {O.\btxfnamespaceshort M. \btxlastnamefont
  {Baksalary}}\btxauthorcolon\ \btxjtitlefont {\btxifchangecase {Particular
  formulae for the {M}oore-{P}enrose inverse of a columnwise partitioned
  matrix}{Particular formulae for the {M}oore-{P}enrose inverse of a columnwise
  partitioned matrix}}.
\newblock \btxjournalfont {Linear Algebra Appl.}, 421(1):16--23,
  2007\ifbtxprintISSN {, \mbox{\btxISSN~\btxISSNfont {0024-3795}}}.
\newblock {\latintext \btxurlfont{https://doi.org/10.1016/j.laa.2006.03.031}}.

\bibitem {MR1927646}
\btxnamefont {\btxlastnamefont {Baksalary},~J.\btxfnamespaceshort K.}
  \btxandshort {.}\ \btxnamefont {G.\btxfnamespaceshort P.\btxfnamespaceshort
  H. \btxlastnamefont {Styan}}\btxauthorcolon\ \btxtitlefont {\btxifchangecase
  {Generalized inverses of partitioned matrices in {B}anachiewicz-{S}chur
  form}{Generalized inverses of partitioned matrices in {B}anachiewicz-{S}chur
  form}}.
\newblock \Btxvolumeshort {.}\ \btxvolumefont {354}, \btxpublisherfont
  {\btxpagesshort {.}\ 41--47}. 2002.
\newblock {\latintext
  \btxurlfont{https://doi.org/10.1016/S0024-3795(02)00334-8}}, Ninth special
  issue on linear algebra and statistics.

\bibitem {MR1987382}
\btxnamefont {\btxlastnamefont {Ben-Israel},~A.} \btxandshort {.}\ \btxnamefont
  {T.\btxfnamespaceshort N.\btxfnamespaceshort E. \btxlastnamefont
  {Greville}}\btxauthorcolon\ \btxtitlefont {Generalized inverses},
  \btxvolumeshort {.}~\btxvolumefont {15} \btxofseriesshort {.}\ \btxtitlefont
  {CMS Books in Mathematics/Ouvrages de Math\'{e}matiques de la SMC}.
\newblock \btxpublisherfont {Springer-Verlag, New York}, \btxeditionnumshort
  {second}{.}, 2003\ifbtxprintISBN {, \mbox{\btxISBN~\btxISBNfont
  {0-387-00293-6}}}.
\newblock Theory and applications.

\bibitem {MR3314341}
\btxnamefont {\btxlastnamefont {Castro-Gonz\'{a}lez},~N.}, \btxnamefont
  {M.\btxfnamespaceshort F. \btxlastnamefont
  {Mart\'{\i}nez-Serrano}}\btxandcomma {} \btxandshort {.}\ \btxnamefont
  {J.~\btxlastnamefont {Robles}}\btxauthorcolon\ \btxjtitlefont
  {\btxifchangecase {Expressions for the {M}oore-{P}enrose inverse of block
  matrices involving the {S}chur complement}{Expressions for the
  {M}oore-{P}enrose inverse of block matrices involving the {S}chur
  complement}}.
\newblock \btxjournalfont {Linear Algebra Appl.}, 471:353--368,
  2015\ifbtxprintISSN {, \mbox{\btxISSN~\btxISSNfont {0024-3795}}}.
\newblock {\latintext \btxurlfont{https://doi.org/10.1016/j.laa.2015.01.003}}.

\bibitem {MR180572}
\btxnamefont {\btxlastnamefont {Cline},~R.\btxfnamespaceshort
  E.}\btxauthorcolon\ \btxjtitlefont {\btxifchangecase {Representations for the
  generalized inverse of sums of matrices}{Representations for the generalized
  inverse of sums of matrices}}.
\newblock \btxjournalfont {J. Soc. Indust. Appl. Math. Ser. B Numer. Anal.},
  2:99--114, 1965\ifbtxprintISSN {, \mbox{\btxISSN~\btxISSNfont {0887-459X}}}.

\bibitem {MR2572505}
\btxnamefont {\btxlastnamefont {Deng},~C.\btxfnamespaceshort Y.} \btxandshort
  {.}\ \btxnamefont {H.\btxfnamespaceshort K. \btxlastnamefont
  {Du}}\btxauthorcolon\ \btxjtitlefont {\btxifchangecase {Representations of
  the {M}oore-{P}enrose inverse of {$2\times 2$} block operator valued
  matrices}{Representations of the {M}oore-{P}enrose inverse of {$2\times 2$}
  block operator valued matrices}}.
\newblock \btxjournalfont {J. Korean Math. Soc.}, 46(6):1139--1150,
  2009\ifbtxprintISSN {, \mbox{\btxISSN~\btxISSNfont {0304-9914}}}.
\newblock {\latintext
  \btxurlfont{https://doi.org/10.4134/JKMS.2009.46.6.1139}}.

\bibitem {MR2641519}
\btxnamefont {\btxlastnamefont {Deng},~C.\btxfnamespaceshort Y.} \btxandshort
  {.}\ \btxnamefont {H.\btxfnamespaceshort K. \btxlastnamefont
  {Du}}\btxauthorcolon\ \btxjtitlefont {\btxifchangecase {Representations of
  the {M}oore-{P}enrose inverse for a class of 2-by-2 block operator valued
  partial matrices}{Representations of the {M}oore-{P}enrose inverse for a
  class of 2-by-2 block operator valued partial matrices}}.
\newblock \btxjournalfont {Linear Multilinear Algebra}, 58(1-2):15--26,
  2010\ifbtxprintISSN {, \mbox{\btxISSN~\btxISSNfont {0308-1087}}}.
\newblock {\latintext \btxurlfont{https://doi.org/10.1080/03081080801980457}}.

\bibitem {MR347845}
\btxnamefont {\btxlastnamefont {Greville},~T.\btxfnamespaceshort
  N.\btxfnamespaceshort E.}\btxauthorcolon\ \btxjtitlefont {\btxifchangecase
  {Solutions of the matrix equation {$XAX=X$}, and relations between oblique
  and orthogonal projectors}{Solutions of the matrix equation {$XAX=X$}, and
  relations between oblique and orthogonal projectors}}.
\newblock \btxjournalfont {SIAM J. Appl. Math.}, 26:828--832,
  1974\ifbtxprintISSN {, \mbox{\btxISSN~\btxISSNfont {0036-1399}}}.
\newblock {\latintext \btxurlfont{https://doi.org/10.1137/0126074}}.

\bibitem {MR1799987}
\btxnamefont {\btxlastnamefont {Gro\ss},~J.}\btxauthorcolon\ \btxtitlefont
  {\btxifchangecase {The {M}oore-{P}enrose inverse of a partitioned nonnegative
  definite matrix}{The {M}oore-{P}enrose inverse of a partitioned nonnegative
  definite matrix}}.
\newblock \Btxvolumeshort {.}\ \btxvolumefont {321}, \btxpublisherfont
  {\btxpagesshort {.}\ 113--121}. 2000.
\newblock {\latintext
  \btxurlfont{https://doi.org/10.1016/S0024-3795(99)00073-7}}, Linear algebra
  and statistics (Fort Lauderdale, FL, 1998).

\bibitem {MR401798}
\btxnamefont {\btxlastnamefont {Hartwig},~R.\btxfnamespaceshort
  E.}\btxauthorcolon\ \btxjtitlefont {\btxifchangecase {Rank factorization and
  {M}oore-{P}enrose inversion}{Rank factorization and {M}oore-{P}enrose
  inversion}}.
\newblock \btxjournalfont {Indust. Math.}, 26(1):49--63, 1976\ifbtxprintISSN {,
  \mbox{\btxISSN~\btxISSNfont {0019-8528}}}.

\bibitem {MR2296072}
\btxnamefont {\btxlastnamefont {He},~C.\btxfnamespaceshort N.}\btxauthorcolon\
  \btxjtitlefont {\btxifchangecase {General forms for {M}oore-{P}enrose
  inverses of matrices by block permutation}{General forms for
  {M}oore-{P}enrose inverses of matrices by block permutation}}.
\newblock \btxjournalfont {J. Nat. Sci. Hunan Norm. Univ.}, 29(4):1--5,
  2006\ifbtxprintISSN {, \mbox{\btxISSN~\btxISSNfont {2096-5281}}}.

\bibitem {MR369383}
\btxnamefont {\btxlastnamefont {Hung},~C.\btxfnamespaceshort H.} \btxandshort
  {.}\ \btxnamefont {T.\btxfnamespaceshort L. \btxlastnamefont
  {Markham}}\btxauthorcolon\ \btxjtitlefont {\btxifchangecase {The
  {M}oore-{P}enrose inverse of a partitioned matrix {$M=(_{B}^{A}$}
  {${}_{C}^{D})$}}{The {M}oore-{P}enrose inverse of a partitioned matrix
  {$M=(_{B}^{A}$} {${}_{C}^{D})$}}}.
\newblock \btxjournalfont {Linear Algebra Appl.}, 11:73--86,
  1975\ifbtxprintISSN {, \mbox{\btxISSN~\btxISSNfont {0024-3795}}}.
\newblock {\latintext
  \btxurlfont{https://doi.org/10.1016/0024-3795(75)90118-4}}.

\bibitem {MR1102137}
\btxnamefont {\btxlastnamefont {Miao},~J.\btxfnamespaceshort
  M.}\btxauthorcolon\ \btxjtitlefont {\btxifchangecase {General expressions for
  the {M}oore-{P}enrose inverse of a {$2\times 2$} block matrix}{General
  expressions for the {M}oore-{P}enrose inverse of a {$2\times 2$} block
  matrix}}.
\newblock \btxjournalfont {Linear Algebra Appl.}, 151:1--15,
  1991\ifbtxprintISSN {, \mbox{\btxISSN~\btxISSNfont {0024-3795}}}.
\newblock {\latintext
  \btxurlfont{https://doi.org/10.1016/0024-3795(91)90351-V}}.

\bibitem {MR3853483}
\btxnamefont {\btxlastnamefont {Mihailovi\'{c}},~B.}, \btxnamefont
  {V.~\btxlastnamefont {Miler~Jerkovi\'{c}}}\btxandcomma {} \btxandshort {.}\
  \btxnamefont {B.~\btxlastnamefont {Male\v{s}evi\'{c}}}\btxauthorcolon\
  \btxjtitlefont {\btxifchangecase {Solving fuzzy linear systems using a block
  representation of generalized inverses: the {M}oore-{P}enrose
  inverse}{Solving fuzzy linear systems using a block representation of
  generalized inverses: the {M}oore-{P}enrose inverse}}.
\newblock \btxjournalfont {Fuzzy Sets and Systems}, 353:44--65,
  2018\ifbtxprintISSN {, \mbox{\btxISSN~\btxISSNfont {0165-0114}}}.
\newblock {\latintext \btxurlfont{https://doi.org/10.1016/j.fss.2017.11.007}}.

\bibitem {MR1488625}
\btxnamefont {\btxlastnamefont {Milovanovi\'{c}},~G.\btxfnamespaceshort V.}
  \btxandshort {.}\ \btxnamefont {P.\btxfnamespaceshort S. \btxlastnamefont
  {Stanimirovi\'{c}}}\btxauthorcolon\ \btxjtitlefont {\btxifchangecase {On
  {M}oore-{P}enrose inverse of block matrices and full-rank factorization}{On
  {M}oore-{P}enrose inverse of block matrices and full-rank factorization}}.
\newblock \btxjournalfont {Publ. Inst. Math. (Beograd) (N.S.)}, 62(76):26--40,
  1997\ifbtxprintISSN {, \mbox{\btxISSN~\btxISSNfont {0350-1302}}}.

\bibitem {MR69793}
\btxnamefont {\btxlastnamefont {Penrose},~R.}\btxauthorcolon\ \btxjtitlefont
  {\btxifchangecase {A generalized inverse for matrices}{A generalized inverse
  for matrices}}.
\newblock \btxjournalfont {Proc. Cambridge Philos. Soc.}, 51:406--413,
  1955\ifbtxprintISSN {, \mbox{\btxISSN~\btxISSNfont {0008-1981}}}.

\bibitem {MR190161}
\btxnamefont {\btxlastnamefont {Rohde},~C.\btxfnamespaceshort
  A.}\btxauthorcolon\ \btxjtitlefont {\btxifchangecase {Generalized inverses of
  partitioned matrices}{Generalized inverses of partitioned matrices}}.
\newblock \btxjournalfont {J. Soc. Indust. Appl. Math.}, 13:1033--1035,
  1965\ifbtxprintISSN {, \mbox{\btxISSN~\btxISSNfont {0368-4245}}}.

\bibitem {MR469933}
\btxnamefont {\btxlastnamefont {Smith},~R.\btxfnamespaceshort
  L.}\btxauthorcolon\ \btxjtitlefont {\btxifchangecase {Moore-{P}enrose
  inverses of block circulant and block {$k$}-circulant
  matrices}{Moore-{P}enrose inverses of block circulant and block
  {$k$}-circulant matrices}}.
\newblock \btxjournalfont {Linear Algebra Appl.}, 16(3):237--245,
  1977\ifbtxprintISSN {, \mbox{\btxISSN~\btxISSNfont {0024-3795}}}.
\newblock {\latintext
  \btxurlfont{https://doi.org/10.1016/0024-3795(77)90007-6}}.

\bibitem {MR2599074}
\btxnamefont {\btxlastnamefont {Tang},~S.} \btxandshort {.}\ \btxnamefont
  {H.\btxfnamespaceshort Z. \btxlastnamefont {Wu}}\btxauthorcolon\
  \btxjtitlefont {\btxifchangecase {The {M}oore-{P}enrose inverse and the
  weighted {D}razin inverse of block {$k$}-circulant matrices}{The
  {M}oore-{P}enrose inverse and the weighted {D}razin inverse of block
  {$k$}-circulant matrices}}.
\newblock \btxjournalfont {J. Hefei Univ. Technol. Nat. Sci.},
  32(9):1442--1444, 1448, 2009\ifbtxprintISSN {, \mbox{\btxISSN~\btxISSNfont
  {1003-5060}}}.

\bibitem {MR1657194}
\btxnamefont {\btxlastnamefont {Tian},~Y.}\btxauthorcolon\ \btxjtitlefont
  {\btxifchangecase {The {M}oore-{P}enrose inverses of {$m\times n$} block
  matrices and their applications}{The {M}oore-{P}enrose inverses of {$m\times
  n$} block matrices and their applications}}.
\newblock \btxjournalfont {Linear Algebra Appl.}, 283(1-3):35--60,
  1998\ifbtxprintISSN {, \mbox{\btxISSN~\btxISSNfont {0024-3795}}}.
\newblock {\latintext
  \btxurlfont{https://doi.org/10.1016/S0024-3795(98)10049-6}}.

\bibitem {MR2490704}
\btxnamefont {\btxlastnamefont {Tian},~Y.} \btxandshort {.}\ \btxnamefont
  {Y.~\btxlastnamefont {Takane}}\btxauthorcolon\ \btxjtitlefont
  {\btxifchangecase {More on generalized inverses of partitioned matrices with
  {B}anachiewicz-{S}chur forms}{More on generalized inverses of partitioned
  matrices with {B}anachiewicz-{S}chur forms}}.
\newblock \btxjournalfont {Linear Algebra Appl.}, 430(5-6):1641--1655,
  2009\ifbtxprintISSN {, \mbox{\btxISSN~\btxISSNfont {0024-3795}}}.
\newblock {\latintext \btxurlfont{https://doi.org/10.1016/j.laa.2008.06.007}}.

\bibitem {MR227186}
\btxnamefont {\btxlastnamefont {Urquhart},~N.\btxfnamespaceshort
  S.}\btxauthorcolon\ \btxjtitlefont {\btxifchangecase {Computation of
  generalized inverse matrices which satisfy specified conditions}{Computation
  of generalized inverse matrices which satisfy specified conditions}}.
\newblock \btxjournalfont {SIAM Rev.}, 10:216--218, 1968\ifbtxprintISSN {,
  \mbox{\btxISSN~\btxISSNfont {0036-1445}}}.
\newblock {\latintext \btxurlfont{https://doi.org/10.1137/1010035}}.

\bibitem {MR566954}
\btxnamefont {\btxlastnamefont {Weidmann},~J.}\btxauthorcolon\ \btxtitlefont
  {Linear operators in {H}ilbert spaces}, \btxvolumeshort {.}~\btxvolumefont
  {68} \btxofseriesshort {.}\ \btxtitlefont {Graduate Texts in Mathematics}.
\newblock \btxpublisherfont {Springer-Verlag, New York-Berlin},
  1980\ifbtxprintISBN {, \mbox{\btxISBN~\btxISBNfont {0-387-90427-1}}}.
\newblock Translated from the German by Joseph Sz\"{u}cs.

\bibitem {MR2517847}
\btxnamefont {\btxlastnamefont {Xu},~Q.}\btxauthorcolon\ \btxjtitlefont
  {\btxifchangecase {Moore-{P}enrose inverses of partitioned adjointable
  operators on {H}ilbert {$C^*$}-modules}{Moore-{P}enrose inverses of
  partitioned adjointable operators on {H}ilbert {$C^*$}-modules}}.
\newblock \btxjournalfont {Linear Algebra Appl.}, 430(11-12):2929--2942,
  2009\ifbtxprintISSN {, \mbox{\btxISSN~\btxISSNfont {0024-3795}}}.
\newblock {\latintext \btxurlfont{https://doi.org/10.1016/j.laa.2009.01.003}}.

\bibitem {MR2993360}
\btxnamefont {\btxlastnamefont {Xu},~Q.}, \btxnamefont {Y.~\btxlastnamefont
  {Chen}}\btxandcomma {} \btxandshort {.}\ \btxnamefont {C.~\btxlastnamefont
  {Song}}\btxauthorcolon\ \btxjtitlefont {\btxifchangecase {Representations for
  weighted {M}oore-{P}enrose inverses of partitioned adjointable
  operators}{Representations for weighted {M}oore-{P}enrose inverses of
  partitioned adjointable operators}}.
\newblock \btxjournalfont {Linear Algebra Appl.}, 438(1):10--30,
  2013\ifbtxprintISSN {, \mbox{\btxISSN~\btxISSNfont {0024-3795}}}.
\newblock {\latintext \btxurlfont{https://doi.org/10.1016/j.laa.2012.08.002}}.

\bibitem {MR2416599}
\btxnamefont {\btxlastnamefont {Xu},~Q.} \btxandshort {.}\ \btxnamefont
  {X.~\btxlastnamefont {Hu}}\btxauthorcolon\ \btxjtitlefont {\btxifchangecase
  {Particular formulae for the {M}oore-{P}enrose inverses of the partitioned
  bounded linear operators}{Particular formulae for the {M}oore-{P}enrose
  inverses of the partitioned bounded linear operators}}.
\newblock \btxjournalfont {Linear Algebra Appl.}, 428(11-12):2941--2946,
  2008\ifbtxprintISSN {, \mbox{\btxISSN~\btxISSNfont {0024-3795}}}.
\newblock {\latintext \btxurlfont{https://doi.org/10.1016/j.laa.2008.01.021}}.

\bibitem {MR3223886}
\btxnamefont {\btxlastnamefont {Yan},~Z.\btxfnamespaceshort Z.}\btxauthorcolon\
  \btxjtitlefont {\btxifchangecase {New representations of the
  {M}oore-{P}enrose inverse of {$2\times 2$} block matrices}{New
  representations of the {M}oore-{P}enrose inverse of {$2\times 2$} block
  matrices}}.
\newblock \btxjournalfont {Linear Algebra Appl.}, 456:3--15,
  2014\ifbtxprintISSN {, \mbox{\btxISSN~\btxISSNfont {0024-3795}}}.
\newblock {\latintext \btxurlfont{https://doi.org/10.1016/j.laa.2012.08.014}}.

\end{thebibliography}
\bibliographystyle{bababbrv}

\vfill\noindent
\begin{minipage}{0.5\textwidth}
 Universit\"at Leipzig\\
 Fakult\"at f\"ur Mathematik und Informatik\\
 PF~10~09~20\\
 D-04009~Leipzig\\
 Germany
\end{minipage}
\begin{minipage}{0.49\textwidth}
 \begin{flushright}
  \texttt{fritzsche@math.uni-leipzig.de\\
   maedler@math.uni-leipzig.de} 
 \end{flushright}
\end{minipage}

\end{document}